\documentclass[11 pt, psamsfonts]{amsart}
\usepackage{amssymb}

\def\ignore#1{\relax}

\def\g{\mathfrak g}
\def\h{\mathfrak h}

\def\R{{\mathbb R}}
\def\Z{{\mathbb Z}}

\def\Q{{\mathbb Q}}
\def\C{{\mathbb C}}
\def\F{{\bf F}}
\def\la{\lambda}
\def\La{\Lambda}

\def\T{{\bf T}}

\def\S{\mathcal S}

\def\N{\mathbb N}
\def\A{\mathcal A}

\def\I{\mathcal I}

\def\P{\mathcal P}

\def\A{\mathcal A}

\def\T{{\bf T}}
\def\a{{\bf a}}
\def\b{{\bf b}}

\def\B{{\mathcal B}}

\def\U{{\bf U}}

\def\ni{\noindent}

\def\v{\vskip 2.5mm}

\def\ignore#1{\relax}

\def\om{\omega}

\def\1{{\bf 1}}

\def\D{{(D)}}

\def\End{{\rm End}}
\def\Hom{{\rm Hom}}
\def\ve{\varepsilon}

\calclayout

\renewenvironment{proof}[1][\proofname]{\par
  \normalfont
  \topsep6\p@\@plus6\p@ \trivlist
  \item[\hskip\labelsep\bfseries
    #1\@addpunct{.}]\ignorespaces
}{%
  \qed\endtrivlist
}

\def\th@plain{%
  \let\thmhead\thmhead@plain \let\swappedhead\swappedhead@plain
  \thm@preskip.5\baselineskip\@plus.2\baselineskip
                                    \@minus.2\baselineskip
  \thm@postskip\thm@preskip
  \itshape
\renewcommand{\labelenumi}{{(\alph{enumi})\quad}}
                        \renewcommand{\labelenumii}{{(\roman{enumii})\ }}
}
\def\th@definition{%
  \let\thmhead\thmhead@plain \let\swappedhead\swappedhead@plain
  \thm@preskip.5\baselineskip\@plus.2\baselineskip
                                    \@minus.2\baselineskip
  \thm@postskip\thm@preskip
  \upshape
}
\def\th@remark{%
  \thm@headfont{\itshape}
  \let\thmhead\thmhead@plain \let\swappedhead\swappedhead@plain
  \thm@preskip.5\baselineskip\@plus.2\baselineskip
                                    \@minus.2\baselineskip
  \thm@postskip\thm@preskip
  \upshape
}

{\theoremstyle{plain}
\newtheorem{theorem}{Theorem}[section]
}

{\theoremstyle{plain}
\newtheorem{prop}[theorem]{Proposition}
}

{\theoremstyle{plain}
\newtheorem{cor}[theorem]{Corollary}
}

{\theoremstyle{plain}
\newtheorem{lemma}[theorem]{Lemma}
}

{\theoremstyle{plain}

}

{\theoremstyle{definition}

}

{\theoremstyle{remark}
\newtheorem{remark}[theorem]{Remark}
}

{\theoremstyle{remark}

}

\numberwithin{equation}{section}

\renewcommand{\labelenumi}{{ \theenumi.}}
\renewcommand{\labelenumii}{{(\alph{enumii})}}

\def\ni{\noindent}

\def\v{\vskip 2.5mm}

\def\la{\lambda}
\def\al{\alpha}

\def\choose #1 #2{\begin{pmatrix}#1\\#2\end{pmatrix}}

\ignore{
\input BoxedEPS
\SetRokickiEPSFSpecial
\HideDisplacementBoxes
\SetEPSFDirectory{./graphics/} }

\ignore{
\input BoxedEPS
\SetTexturesEPSFSpecial
\HideDisplacementBoxes
\SetEPSFDirectory{:graphics:}
}

\usepackage{amsfonts}
\usepackage{amsmath}
\usepackage{graphics}
\begin{document}
\def\e{\epsilon}
\def\ga{\gamma}
\def\al{\alpha}
\def\g{\mathfrak g}
\def\a{\alpha}
\def\om{\omega}
\def\b{\beta}
\def\e{\epsilon}
\def\R{\mathbb{R}}
\def\ga{\gamma}
\def\so{{\mathfrak s}{\mathfrak o}}
\def\sp{{\mathfrak s}{\mathfrak p}}
\def\la{\lambda}
\def\Ga{\Gamma}
\def\DN{D(N)}
\def\l{{\bf l}}
\def\ve{\varepsilon}
\def\diag{{\rm diag}}
\def\h{\mathfrak h}
\def\ZZ{\mathbb{Z}}
\def\N{\mathbb{N}}
\def\Q{\mathbb{Q}}
\def\CC{\mathbb{C}}
\def\A{\mathcal{A}}
\def\P{\mathcal{P}}
\def\S{\mathcal{S}}
\def\C{\mathcal{C}}
\def\D{\mathcal{D}}
\def\Cc{\hat{\mathcal{C}}}
\def\B{\mathcal{B}}
\def\Grin{Gr(O(\infty))}
\def\CSinf{\cup \CC S_n}
\def\Cinf{\cup \C_n}
\def\CN{\cup \bar\C_n(N)}
\def\I{{\mathcal I}}
\def\IONN{{\mathcal IO}_N}
\def\INN{{\mathcal I}_N}
\def\ISpNN{{\mathcal ISp}_N}
\def\ION{{\mathcal IO}_N}
\def\IN{{\mathcal I}_N}
\def\ISpN{{\mathcal ISp}_N}
\def\UU{{\bf U}}
\def\ni{\noindent}
\def\U{U_q\g}
\def\La{\Lambda}
\def\dim{\mbox{dim}}
\def\v{\vskip .2cm}
\def\End{\mbox{End}}
\def\Hom{\mbox{Hom}}
\def\Vf{V_{m\e}\otimes V^{\otimes f}}
\def\Wm{{W^{(m)}}}
\def\Wl{{W^{(l)}}}
\def\x{{\bf x}}
\def\y{{\bf y}}
\def\j{{\bf j}}
\def\bCcn{\bar \C_n}
\def\bCc{\bar \C}
\def\sh{{\rm sh}}
\def\T{{\mathcal T}}
\def\F{{\mathcal F}}
\newtheorem{observation}[theorem]{Observation}
\newtheorem{definition}[theorem]{Definition}

\title{\bf{Quotients of Representation Rings}}

\author
{Hans Wenzl}
\thanks{Partially supported by NSF grant DMS 0302437}

\begin{abstract}
We give a proof, using so-called fusion rings and $q$-deformations of
Brauer algebras
that the representation ring of an orthogonal or symplectic
group can be obtained as a quotient of a  ring
$\Grin$. This is obtained here as a limiting case for analogous
quotient maps for fusion categories, with the level going to $\infty$.
This in turn allows a detailed description of the quotient map
in terms of a reflection group.
As an application, one obtains a general description of the branching
rules for the restriction of representations of $Gl(N)$
to $O(N)$ and $Sp(N)$ as well as detailed information
about the structure of the $q$-Brauer algebras in
the nonsemisimple case for certain specializations.
\end{abstract}

\maketitle

It is well-known that one can study the combinatorics of the
finite dimensional representations of
the general linear groups $Gl(N)$ in a uniform way
essentially independently of the dimension $N$, using the
representation theory of the symmetric groups.
A similar approach is possible for orthogonal and symplectic groups.
In particular, it is possible to obtain their Grothendieck rings
as quotients from a large ring, here denoted by $\Grin$
(see e.g. \cite{KT}). More recently, quotients $Gr(G)_\ell$
of the Grothendieck ring $Gr(G)$
of a semisimple Lie group $G$, depending on a positive integer $\ell$,
 have been studied which have
only finitely many simple objects up to isomorphism. They originally
arose in mathematical physics and are usually called fusion rings.
One of the observations in this paper is that there is more than
just a formal similarity between these two quotients:

We show that the quotient map from $\Grin$ onto the Grothendieck
ring $Gr(Sp(N))$ of a symplectic group can be obtained as
a limit of the quotient map $Gr(Sp(M))$ to $Gr(Sp(M))_{\ell(M)}$
for $M\to\infty$; here $\ell(M)$ depends on $M$ in a simple linear way.
A similar result also holds for orthogonal groups.
This allows, among other things, for a simple explicit description
of the quotient map in terms of a reflection group.
One of the applications is an explicit formula for the restriction
multiplicities for representations from $Gl(N)$ to $O(N)$ and to
$Sp(N)$ respectively; see below for more details and related results.
As another application, we also obtain results about the structure
of a $q$-deformation $\Cc_n(r,q)$ of Brauer's centralizer algebras
$\C_n(x)$ when $r=q^{-N-1}$, with $N>0$ even. In our approach,
we get these results as a limiting case of rather deep results
about tilting modules of quantum groups. This also makes the
appearance of parabolic Kazhdan-Lusztig polynomials fairly natural.

Here is the contents of our paper in more detail. The first section
reviews results about so-called fusion rings, which are certain
quotients of the representation rings of semisimple Lie algebras,
depending on a positive integer $\ell$. This is fairly elementary.
We can deduce from this one of our main results, the already
mentioned explicit description of the quotient maps from $\Grin$
onto $Gr(Sp(N))$ and $Gr(O(N))$ respectively. The only result
which is not proved by elementary methods is the so-called level-rank
duality, see Prop. \ref{fusionsym}. A very simple proof can be given
after expressing the corresponding fusion rings
in terms of $q$-versions of Brauer's centralizer algebras.
Section 2 reviews the necessary results about these algebras.
In Section 3 we define semirings in terms of these algebras and
eventually show these semirings are related to the fusion rings defined
in the first section. Section 4 then proves the
restriction rules from $Gl(N)$ to $O(N)$ and $Sp(N)$, respectively.
Finally, in Section 5, we first review how results by Soergel and Hu
determine the structure of the image of
the algebras $\Cc_n(q^{-N-1},q)$, $N$ even, in a certain tensor space.
Using the already mentioned level-rank duality, we extend this
to the whole algebras in the generic case, i.e. for all but finitely
many values of $q$.

This paper deals with some very familiar and heavily studied subjects.
So, not surprisingly, many related results have appeared before.
We just mention a few: This paper was motivated by the paper \cite{EW}
where restriction rules have been studied by a completely different
approach. After a first version of our paper was finished, we also
learned of the work \cite{KT} with very similar results, and an
approach much closer to the one used here, and other combinatorial
approaches, including \cite{Ki} and \cite{Su}.
Nevertheless, we think
that the connection to fusion rings and tilting modules in our
paper provides additional insight.
In the version of the paper which was first submitted in December 2006,
there was also a conjecture about the structure of Brauer's centralizer
algebra in the nonsemisimple case, inspired by Soergel's work on
tilting modules.
This conjecture was also obtained, independently,
in \cite{CdVM} and, more recently, a proof has been posted in \cite{Mt}.
In addition,  the  paper \cite{Hu} has appeared which proves one of
the assumptions for our approach. We decided to replace the section
dealing with the conjecture with a new section which does prove
the conjecture in a special case for the $q$-deformation of the Brauer
algebras. This should give the same result as the one in \cite{Mt}
if $q\to 1$, but the method would be completely different.
The rest of the paper has the same results as the first
version. Some substantial changes have been made in the hope of making
the presentation more conceptual and clearer, but no new material
has been used for that.

$Acknowledgement:$ I would like to thank the referee for his/her remarks
which led to an improved presentation of the material.

\section{Fusion rings}

\ignore{
\subsection{Multiplicative characters}
We call any ring homomorphism $\psi:\Grin \to \CC$
a {\it multiplicative character} of $\Grin$.
Using the duality with $G=O(N)$ or $G=Sp(N)$,
it is easy to construct examples of characters
for $Gr(\C(N))$: If $g\in G$, $p(N)\in\C_n(N)$ and $Tr$ is the usual
trace on $V^{\otimes n}$,
we define $\chi=\chi_g$ by $\psi([p[N]])=Tr(g^{\otimes n}\Phi(p(N)))$.
It is easy to check that this defines a multiplicative character on
$Gr(\C(N))_\pm$, and hence also on $\Grin_\pm$ (see Prop. \ref{quot}).
\begin{definition}\label{Gcharacters}
Let $G=O(N)$ or $Sp(N)$, and let $g\in G$. We call a character $\psi$ of
$\Grin$ a $G$-character if there exists a $g\in G$ such that $\psi=\psi_g$.
The ideal $\I^G$ is defined as $\I^G=\{ d\in\Grin,\ \psi_g(d)=0$ for all
$g\in G\}$.
\end{definition}}

\subsection{General construction}\label{generalconstruction}
The material here is well-known. For some background, see e.g.
the books \cite{FMS} and \cite{Kac}.
Let $G$ be a compact simple
Lie group. Let $\g$ be its
Lie algebra  with Cartan subalgebra $\h$, and let
$\h^*$ be the linear dual of $\h$. As in \cite{Kac} we denote
the $coroot$ $lattice$ by $M$.
We normalize
the invariant bilinear form on  $\h^*$ such that $(\al, \al)=2$ for a
long root $\al$. Identifying $\h^*$ with $\h$ via this form, we can
assume that the coroot lattice $M$ is a sublattice of the weight lattice
$M^*$. Moreover, the weight lattice is given by
$M^*=\{ \y\in\h^*, (\y,\x)\in\Z$ for all $\x\in M\}$.
For a given dominant weight $\la$, we denote by $V_\la$ a simple
$\g$ module with highest weight $\la$, and by $[V_\la]$ its
equivalence class in the Grothendieck semiring of $\g$.
We will explicitly write down $M$ and $M^*$ for the Lie algebra of
an orthogonal or a symplectic group below.

For any $\y\in\h^*$ we obtain a multiplicative character $e^\y$
of $\h$ via $e^\y(\x)=exp(2\pi i(\y,\x))$.
Fix $\ell\in\N$ and let $\ell^{-1}M^*=\{ \y\in\h, \ell\y\in M^*\}$.
Then $\ell^{-1}M^*/M$ is a finite group whose simple
characters are given by $\{ e^\y, \y\in M^*/\ell M\}$.
Let $W$ be the Weyl group associated to $G$.
We define, as usual,
the Weyl character $\chi^\la(\x)=(\sum_{w\in W} \ve(w)e^{w(\la+\rho)}(\x))/
(\sum_{w\in W} \ve(w)e^{w(\rho)}(\x))$,
where $\rho$ is half the sum of the positive
roots of $G$.
Let $\Wl$ be
the associated affine reflection group. The latter is the semidirect
product of $W$ with $\ell M$, where the latter group acts via
translations $t_{\x}\y=\y+\x$ on $\h^*$. It is well-known that
$\Wl$ has a fundamental domain  $\bar D_\ell$,
the closure of the set $D_\ell$ consisting of elements $\y\in\h^*$ satisfying
$(\al_i,\y)> 0$ and $(\theta,\y)< \ell$; here the $\al_i$ are the
simple roots of $G$, and $\theta$ is the maximum long root of $G$.
It is well-known that $\Wl$ is a reflection group, generated by
the reflection $s_0$ in the hyperplane given by
$(\theta,\y)=\ell$ and the simple reflections of $W$.
The hyperplanes for these generating reflections contain the boundary
simplices of $D_\ell$.
We define a second action of $\Wl$ on $\h^*$ by
$w.\la=w(\la+\rho)-\rho$ for $\la\in\h^*$.

\ignore{
The following lemma is well-known and is easy
to prove; as usual $\ve(w)$ denotes the sign of the element $w$ of a
reflection group.
\begin{lemma}\label{reduction}
$a_W(e^{w(\y)})=\ve(w)\a_W(e^\y)$ for all $w\in \Wl$
and all $\y\in \ell^{-1}M$.
\end{lemma}
}

\begin{prop}\label{fusionring} Let $G$ be a connected and simply
connected simple Lie group.
Let $I_\ell(G)$ be the ideal in the representation
ring of $G$ generated by all those characters which annihilate
$\ell^{-1}M^*$. Then $Gr(G)_\ell=Gr(G)/I_\ell(G)$ has a $\Z$-basis
$([V_\la])$ labeled
by all those dominant weights $\la$ for which $\la+\rho$ is in $D_\ell$.
Moreover, we have $[V_{w.\la}]=\ve(w)[V_\la]$ modulo $I_\ell(G)$
for any $w\in\Wl$ which
maps the dominant weight $\la$ to a dominant weight $w.\la$.
The spectrum of the abelian ring  $Gr(G)_\ell$ is given by the maps
$\phi_\mu: [V_\la]\mapsto \chi^\la((\mu+\rho)/\ell)$,
where $\mu$ goes over all dominant weights of $G$ such that
$\mu+\rho\in D_\ell$.
\end{prop}

$Proof.$ This is well-known.
We give the simple arguments for
the reader's convenience. Define for $\y\in M^*$ the map
$a_W(e^\y)=\sum_{w\in W}\ve(w)e^{w(\y)}$.
If $\la$ is a dominant weight,
$a_W(e^{\la+\rho})$ is the numerator in the Weyl character formula
for the simple  module with highest weight $\la$.
Moreover, for such $\la$,
we can find $\x\in\ell M$ and $w\in W$ such that $wt_\x(\la +\rho)$
is in $\overline{D_\ell}$ (The reader not familiar with this should
have little problems checking this for the relevant cases which
will be listed below).
It follows that $a_W(e^{\la+\rho})=\ve(wt_\x)a_W(e^{wt_\x(\la+\rho)})$;
the latter expression is zero if $wt_\x(\la+\rho)$ is not in $D_\ell$,
i.e. if it lies on one of its boundary simplices.

Finally, it is clear that each $\phi_\mu$ defines a homomorphism
from $Gr(G)_\ell$ into $\CC$. Let $d(\mu)=\chi^\mu(\rho/\ell)$
be the so-called $q$-dimension of $V_\mu$. Then
$s_{\la\mu}=\chi^\la((\mu+\rho)/\ell)d(\mu)$ is an entry of the
$S$-matrix $S=(s_{\la\mu})$ in \cite{KP}. It is shown there that this
matrix is invertible (or see, e.g. \cite {TW2}, Prop. 9.8.1 for the
simple argument). Hence the $\phi_\mu$ are linearly independent
and they exhaust all one-dimensional characters of $Gr(G)_\ell$.

\begin{remark}\label{Grplus} 1. We define $Gr(G)_{\ell, +}$ to be
the $\N$-linear span of the elements $\la\in Gr(G)_\ell$ for which
$\la+\rho$ is in $D_\ell$. We shall see later that this does indeed
 define a semiring.

 2. With some care, one can similarly also define quotients $Gr(G)_\ell$
 if $G$ is not necessarily connected or simply connected. This will be
 done below for orthogonal groups.
\end{remark}

\subsection{Orthogonal case}
Let us first explicitly describe the set-up of the last subsection
for the odd-dimensional orthogonal Lie algebra, i.e. for Lie type $B_m$.
For type $B_m$, we can take
$M= \{ \x\in \Z^m, 2|\sum x_i\}$ and
$M^*=\Z^m\cup (\j+\Z^m)$, where $\j=(1/2,1/2,\ ...,\ 1/2)\in\R^m$. Here
the invariant inner products on $M$ and $M^*$ are given by the
usual scalar product on $\R^m$. Its Weyl group $W(m)$ is isomorphic
to the semidirect product of $(\Z/2)^m$ with the symmetric group $S_m$.
It has a natural action on $\h^*\cong \R^m$
where $S_m$ permutes the coordinates
and $(\Z/2)^m$ acts via possible sign changes on the coordinates.
Its maximum longest root is given by $\theta(\x)=x_1+x_2$ for $\x\in\R^m$.
The corresponding affine Weyl group $\Wl(m)$ has as fundamental domain
the closure of the set
$D_\ell=\{ x\in\R^m, x_1>x_2>\ ...\ x_m>0$ with $x_1+x_2<\ell\}$.
Finally, we set $\rho(m)$ to be equal to half of the sum of all positive
roots; in our notation, that means $\rho(m)=(m+\frac{1}{2}-i)_i\in\R^m$.

For the Lie algebra $so_{2m}$, we have the same lattices $M$ and $M^*$
and the same longest root $\theta$ as for Lie type $B_m$. The Weyl group
$W(m)$ now is isomorphic to the semidirect product $(\Z/2)^{m-1}\ltimes S_m$,
where $(\Z/2)^{m-1}$ is identified with the subgroup of $(\Z/2)^m$
consisting of elements with an even number of nonzero entries.
We define the set $D_\ell$ for the affine Weyl group $W^{(\ell)}(m)$
as for type $B_m$ except for the inequality $...\ >x_{m-1}>|x_m|$;
its closure is again a fundamental domain.

\subsection{The semiring $\Grin_+$}\label{Grin}
Let $M$ and $M^*$ be as in the last
subsection, for the odd-dimensional case $SO(2m+1)$.
Only the weights in $M^*$ with integer coefficients are also weights
of the corresponding orthogonal group. These weights are invariant
under the action of the group $\Wl(m)$, hence we obtain a well-defined
subgroup of the quotient $Gr(Spin(2m+1))_\ell$, as defined in Prop.
\ref{fusionring}. In particular,
the dominant weights are given by vectors with nonnegative nonincreasing
integer coefficients. We shall often identify such a dominant weight $\la$
with a Young diagram, denoted by the same symbol, which has $\la_i$ boxes
in its $i$-th row. This allows an obvious injection of the dominant
weights of $SO(2m+1)$ into the dominant weights of $SO(2(m+1)+1)$.
It turns out that for given dominant weights $\la$ and $\mu$
these embeddings are compatible with the tensor product
rules for large enough $m$. More precisely, we have the following
well-known proposition (see e.g. \cite{KT}; we will give another proof
in Lemma \ref{genericring}).

\begin{prop}\label{stablerules}
Let $V_\la,\ V_\mu$ be simple $SO(2m+1)$ modules labeled by
Young diagrams $\la$ and $\mu$, and let
$V_\la\otimes V_\mu=\bigoplus d_{\la\mu}^\nu(m) V_\nu$
be the decompostion of the tensor product into a direct sum of simple
modules. Then $ d_{\la\mu}^\nu(m) =  d_{\la\mu}^\nu$ independently of
$m$ for $m$ sufficiently large.
In particular, we can define a semiring $\Grin_+$ with an $\N$ basis
$\{[\la] \}$ labeled by Young diagrams whose multiplication is given
by the structure coefficients $d_{\la\mu}^\nu$.
\end{prop}

\subsection{Limiting reflection group}\label{lrefl}
With the notations of last section for type $B_m$,
let $\ell=N+2m-1$ for $N$ a fixed
positive integer. It will turn out that we can construct a limiting
Grothendieck ring for $m\to\infty$, similarly as we defined
the ring $\Grin$ in the last section. It will turn out that for
given $\la$ and large $m$, only certain reflections of $\Wl(m)$
will be relevant to map $\la$ into $D_\ell$.
Using only these reflections, with the standard embedding
$W(m)\subset W(m+1)$, we obtain a limiting reflection
group $W=W{(\infty)}$.
This group $W$ has Coxeter type $D_\infty$.
It is generated by the group $S_\infty$ of finite
permutations on a countable set together with one additional reflection
$s_o$. It can be conveniently described by
a faithful representation on the space of
sequences, on which $S_\infty$ acts via permutations of the entries
of the sequence and where the reflection $s_0$ acts by
$s_0(l_1,l_2,l_3,\ ...)=(-l_2,-l_1,l_3,\ ...)$.
In particular, we can define the sign $\ve(w)$
of an element $w\in W$ to be $(-1)^k$, where $k$ is the number of factors
if we write $w$ as a product of simple reflections.
Let $\rho$ be the sequence $(1-N/2-i)_{i\in\N}$. Then we define
a second action on sequences $\la$ as before by
$w.\la=w(\la+\rho)-\rho$.

We also define $\DN$ to be the set of all Young diagrams with
$\leq N$ boxes in the first two rows. Similarly as $Gr(G)$ can be
obtained as a quotient of $\Grin$, we can now define a ring
which contains all the orthogonal fusion rings as quotients.
The algorithm in the proof below already appeared before, see
\cite{Su}.

\begin{lemma} \label{redlem}  Let $N$ be a fixed positive integer,
and let $\ell=N+2m-1$, with also $m$ being a positive integer.
Let $\la\in\La$ be a Young diagram with $\leq m$ rows.

(a) There exists an element
$w\in W$ such that either $w.\la\in\DN$ or
$\la+\rho$ is fixed by some reflection in $W$.
Moreover, if $w.\la\in\DN$,
it is contained in $\la$, i.e. $(w.\la)_i\leq \la_i$ for all $i$.

(b) If $\la_1+\la_2<2(N+m+1)$ and if $\tilde w\in \Wl{(m)}$
is such that $\tilde w.\la\in D_N$,
then $\tilde w.\la =w.\la$, with $w$ as in (a); in particular,
$\tilde w.\la$ does not depend on $m$.
\end{lemma}

$Proof.$  Let us first consider the situation in part (b).
Set $\x=\la+\rho(m)$. There is nothing to show if $x_1+x_2\leq \ell$.
By assumption, we have $x_1+x_2<2\ell$. Hence, if $\y=s_0(\x)$,
we have $0<\y_1+\y_2<\x_1+\x_2$. If $\y$ does not have
two identical coordinates (which does imply that  $\x$ is fixed
by a reflection), we can apply a permutation $w$ to $\y$
such that $w(\y)$ has strictly decreasing coordinates, i.e.
$w(\y)-\rho=\tilde\la$ is
a Young diagram with $\tilde\la_i\leq\la_i$ for all $i$ and
$\tilde\la_1+\tilde\la_2<\la_1+\la_2$.
Iterating this process, we will end up with a diagram as
described in part (a).

To prove the remaining parts of this lemma, it suffices to check
that $s_0.\la= \la-(\la_1+\la_2-N-1)(1,1,0, ...)$ both for sufficiently
large $m$ as well as for $W=W(\infty)$. This is straightforward.
Hence the
algorithm in the last paragraph works as well in the setting of
$W=W(\infty)$, and yields the same result. $\Box$

\begin{prop}\label{orthogonalreduction} For each $N\in \N$
there exists an ideal $\IN=\I_{O(N)}$ such that the
quotient $\Grin/\IN$ of $\Grin$ has a $\Z$-basis labelled
by the set $\DN$ of Young diagrams $\la$ satisfying $\la_1+\la_2\leq N$.
For other Young diagrams $\mu$, we either have $[\mu]\in\IN$, or
there exists a $w\in W$ such that $w.\mu\in \DN$ and
$[\mu]\cong \ve(w)[w.\mu]$ mod $\IN$.
\end{prop}

$Proof.$ We have seen that the product $[\la][\mu]$ of elements
$[\la],[\mu]\in\Grin_+$ is given by the decomposition of the
tensor product $V_\la\otimes V_\mu$ of irreducible $SO(2m+1)$-modules,
for $m$ sufficiently large. It only remains to show that the reduction
modulo $\IN(m)$ does not depend on $m$ provided $m$ is sufficiently large.
But it follows from the algorithm in the proof of
Lemma \ref{redlem} that the fact whether $[\la]\in \IN(m)$ and,
if not, the value of $w.\la\in\DN$ is independent of $m$, provided
it is large enough. $\Box$

\subsection{Symplectic case}\label{sympsect}
In the symplectic case, Lie type $C_m$, we can take $M=2\Z^m$ and
$M^*=\Z^m$, with the inner product given by {\rm one half} of the usual
scalar product on $\R^m$. The Weyl group $W(m)$ is the same as for type $B_m$,
while the maximum long root $\theta$ now is given by $\theta(\x)=x_1$.
We again define the affine reflection group $\Wl(m)$ as the semidirect
product of $\ell M$ with $W$. It now has a fundamental domain given
by the closure of the set
$D_\ell=\{ x\in\R^m, x_1>x_2>\ ...\ x_m>0$ with $x_1<\ell\}$.
In this setting, half the sum of the positive roots is equal to
$\rho(m)=(m+1-i)_i\in\R^m$.

We will again define a limiting reflection group $W=W(\infty)$,
this time of type $B_\infty$, which is generated by
$S_\infty$ and an additional reflection $s_0$. It has a
faithful action on the set of sequences given by
$s_o(l_1,l_2,\ ...)=(-l_1,l_2,\ ...)$ and by the obvious
action of $S_\infty$. We also define $\rho$ to be the
sequence $\rho=(N/2-i)_i$.

\begin{prop}\label{sympreduction} For each even $N\in \N$
there exists an ideal $\IN=\I_{Sp(N)}$ such that
the quotient $\Grin/\IN$ of $\Grin$ has a $\Z$-basis labelled
by the set $\DN=D(Sp(N))$ of Young diagrams $\la$ satisfying $\la_1\leq N/2$.
For other Young diagrams $\mu$, we either have $[\mu]\equiv 0$ mod $\IN$, or
there exists a $w\in W$ such that $w.\mu\in \DN$ and
$[\mu]\cong \ve(w)[w.\mu]$ mod $\IN$.
\end{prop}

$Proof.$ The proof follows exactly the same pattern as
the one for Prop. \ref{orthogonalreduction}.
So we only give an outline of the proof. One defines ideals
$\IN(m)$ for the  symplectic groups $Sp(2m)$ with $\ell=N/2+m+1$
as in Prop. \ref{fusionring}.
One then proves the analog of Lemma \ref{redlem} in this
setting.

\subsection{Full orthogonal groups}\label{fullortho}
We will need to extend the discussion
above for connected groups to the full orthogonal group $O(M)$, which
is a semidirect product of $SO(M)$ with $\Z/2$. Its simple modules are
labeled by Young diagrams with at most $M$ boxes in the first two columns;
using the notation $\la_i'$ for the number of boxes in the $i$-th column,
we can formulate this condition by $\la_1'+\la_2'\leq M$.
We define the map $t$ on Young diagrams by $t(\la)'_1=M-\la_1'$ and
$t(\la)_i'=\la_i'$ for $i>1$. Observe that $t$ defines a permutation
of the irreducible representations of $O(M)$ of order two.
We can now easily describe the restriction rule from $O(M)$ to $SO(M)$
as follows:

$(i)$ If $t(\la)\neq \la$, the irreducible $O(M)$ modules $V_\la$
and $V_{t(\la)}$ labeled
by $\la$ and $t(\la)$ respectively are isomorpic irreducible
$SO(M)$-modules. One of $\la$ or $t(\la)$ will have less than
$M/2$ rows and will give us the highest weight of the $SO(M)$-module
$V_\la$.

$(ii)$ If $t(\la)=\la$, then $M$ is even and $\la$ has exactly $M/2$ rows.
In this case the irreducible $O(M)$ module $V_\la$ decomposes into
the direct sum of two irreducible $SO(M)$ modules with highest weights
$\la$ and $\tilde\la$, where $\tilde\la$
coincides with $\la$ except for the $m$-th coordinate, which is
equal to $-\la_m$.

We can now extend the quotient map from $Gr(SO(M))$ to $Gr(O(M))$
as follows:

\begin{lemma}\label{fusionsymm}
The ring $Gr(O(M))$ has a quotient $Gr(O(M))_\ell$  with a $\Z$ basis given by
Young diagrams $\la$
with $\leq M$ boxes in the first two columns, and with $\leq \ell + 2-M$
boxes in the first two rows.
This extends the quotients constructed above
for $Gr(SO(M))$, with the same restriction rules.
In particular, the quotient map is determined by
the same group $W$ as in the $SO(M)$ case, which now acts only
on the first $m$ rows.
\end{lemma}

$Proof.$ If $M$ is odd, $O(M)\cong SO(M)\times \Z/2$. The action of
the generator of $\Z/2$ on $V_\la$ is given by the scalar $(-1)^{|\la |}$.
It is easy to see that this is compatible with the action of the group
$W$. The case is more complicated for $M=2m$ even.
If $\la$ has less than $m$ rows, we can apply $W$ as in the odd case.
If $\la$ has exactly $m$ rows, it can happen that $w.\la$ may have
less than $m$ rows. In this case, it follows from the character formula
that $[\la]\cong \ve(w)([w.\la ] +[t(w.\la)])$.
Finally, if $\la$ has more than $m$ rows, we have $[w.\la]=[t(w.(t(\la)))]$.

\begin{prop}\label{fusionsym} (Level-rank duality)
The transpose map $\la\mapsto \la'$, where $\la'$
is the transposed Young diagram obtained from $\la$, induces an
isomorphism between $Gr(O(N))_\ell$ and $Gr(O(\ell+2-N))_\ell$
and between $Gr(Sp(N))_\ell$ and $Gr(Sp(2\ell -2-N))_\ell$.
\end{prop}

This proposition will be shown at the end of Section \ref{Grothendiecksemi};
see a brief discussion about this result at the beginning of that section.
For the following theorem, we define the transposed action of $w\in W$ on a
Young diagram by $w'.\la=(w.\la')'$, i.e. we define the action via
the columns instead of the rows.

\begin{theorem}\label{classicsym}
The transpose map $\la\mapsto \la'$ induces an isomorphism
between the rings $Gr(O(N))$ and $\Grin/\I_{O(N)}$. In particular, $Gr(O(N))$
is a quotient of $\Grin$, with the quotient map explicitly given
by Prop \ref{orthogonalreduction} with respect to the transposed action
of $W$. Similarly, we have isomorphisms between
$Gr(Sp(N))$ and a quotient $\Grin/\I_{Sp(N)}$, using the notations of
Section \ref{sympsect}.
\end{theorem}

$Proof.$ The transpose maps  $\Grin/\I_{O(N)}$ to a quotient of $\Grin$
which has a basis labeled by the same Young diagrams which also
label the irreducible representations of $O(N)$. By definition of
$\Grin$ and  $\Grin/\I_{O(N)}$, the multiplication of two elements
$[\la]$ and $[\mu]$ in the latter quotient is already determined
by the one in $Gr(O(\ell+2-N))_\ell$ for $\ell$ sufficiently large,
see Lemma \ref{redlem} and Prop. \ref{orthogonalreduction}.
By Prop. \ref{fusionsym}, the latter semiring is isomorphic to
$Gr(O(N))_\ell$.
Again, for $\ell$ sufficiently large, the product is given by
the  usual tensor product rules of $O(N)$. The same proof works
as well in the symplectic case.

\section{Brauer algebras}\label{brauersec}
\subsection{Basics}\label{basics} (See \cite{Brauer})
 Let $V=\CC^N$ be the vector representation of
the orthogonal group $O(N)$. Then, using the usual bilinear form
on $V=\CC^N$, we can identify $V^*$ with $V$, and $\End(V^{\otimes n})$
with $(V^*)^{\otimes 2n}\cong (V^{\otimes n})^*$. Moreover, under this
isomorphism, the linear space $\End_{O(N)}(V^{\otimes n})$ of linear
maps commuting with the $O(N)$ action on $V^{\otimes n}$ is isomorphic
to the linear space of $O(N)$-invariant elements in $(V^{\otimes n})^*$.
It was shown by Brauer that this space is spanned by all possible
maps obtained by forming the product of $n$ inner products of the $2n$
tensor factors of an element in $V^{\otimes 2n}$. As an example, for $n=2$,
we have three such maps, which map a vector
$v_1\otimes v_2\otimes v_3\otimes v_4$ to the products
$(v_1,v_3)(v_2,v_4)$, $(v_1,v_4)(v_2,v_3)$ and $(v_1,v_2)(v_3,v_4)$
respectively. In general, these maps are given by a partition
of the set $\{ 1,2,\ ...,\ 2n\}$ into $n$ disjoint subsets $\{ i_r,j_r\}$,
$1\leq r\leq n$ which describe the map
\begin{equation}\label{brauerbasic}
v_1\otimes v_2\otimes\ ...\ \otimes v_{2n}\quad \mapsto \quad \prod_{r=1}^n
(v_{i_r},v_{j_r}).
\end{equation}

\subsection{Definitions}\label{defsec} (See \cite{Brauer},
and also \cite{w2}.) The discussion in the previous section
motivated the definition of an abstract algebra $\C_n$
which can be defined over the ring $\Z[x]$ of polynomials with
integer coefficients; we will usually define it over
the rational functions $\Q(x)$. For each map
in Eq \ref{brauerbasic} we define a graph
with $2n$ vertices and $n$ edges, where  we put the first $n$ vertices in
a lower, and  the remaining $n$ vertices in an upper row
such that the $(n+i)$-th vertex is above the $i$-th vertex for
$1\leq i\leq n$. Now the $r$-th edge connects
the vertices $i_r$ and $j_r$ in our graph, and
and horizontal edges (i.e. edges which connect vertices which are
in the same row) should be drawn slightly concave.
Multiplication of a
graph $a$ with a graph $b$ is given by putting $a$ on top of $b$,
i.e. by identifying the $n$ lower vertices of $a$ with the $n$ upper
vertices of $b$. The composite graph is again a basis graph, except
that there may be cycles, i.e. components of the graph which are
no longer connected with any upper or lower vertex (see example below).
The product $ab$ is then defined to be the
graph obtained by removing all cycles from the composite graph,
multiplied by $x^c$, where $c$ is the number of cycles.
As an instructive example, let
$e$ denote the graph in $\C_2$ with horizontal upper and lower edges.
Then putting two of such graphs on top of each other produces
a cycle (thanks to the concave drawing) between two horizontal edges
and one obtains $e^2=xe$ (see e.g. \cite{w2}, p. 181).

Observe that $\C_n$ contains the group algebra of the symmetric group
$ S_n$ as a subalgebra; it is spanned by graphs which connects
the $i$-th lower vertex to an upper vertex, say the $\pi(i)$-th vertex,
for $i=1,2,\ ...\ n$, where $\pi\in S_n$; here we label both
the lower and the upper vertices by the numbers 1 until $n$.
We will call such graphs
permutation graphs.
It is not hard to see that $\C_n$ is generated by $S_n$ and
$e_1=e\otimes 1_{n-2}$.

Moreover, we  define for graphs $a\in\C_n$, $b\in\C_m$
a new graph $a\otimes b\in\C_{n+m}$
by putting the graph $b$ to the right of the graph $a$. Then $a^{\otimes k}$
is equal to $a\otimes a\otimes\ ...\ \otimes a$ ($k$ times) and
the identity $1_k$ for $\C_k$ is given by the graph with $k$ vertical edges.

Obviously, the definitions above go through if we replace $\Z[x]$ by
any polynomial ring $R[x]$ for a unital abelian ring $R$
or by any quotient field $F(x)$ for $F$ a field. We will only
consider characteristic 0 in this paper.
Similarly, if $N\in R$ for a unital abelian ring,
we define the  algebra $\C_n(N)$ over $R$
via the same graphs by substituting $x=N$ in the multiplication.
We denote by $\C_n(N)$ the $\Q$ algebra spanned by the basis graphs,
with $x=N\in\Q$.

\subsection{Centralizers}\label{centralizers}
We have already seen how a basis graph
is related to a map in Eq \ref{brauerbasic}. Combining this with
the canonical  isomorphism between $(V^*)^{\otimes 2n}$ and
$\End(V^{\otimes n})$, we obtain for each basis graph
a map in $\End_{O(N)}(V^{\otimes n})$. E.g. for $n=2$
the three maps mentioned in Section \ref{brauerbasic}
correspond to the identity map, the permutation map $v_1\otimes v_2
\mapsto v_2\otimes v_1$ and to the map
\begin{equation}\label{Emap}
\Phi(e):\ v\otimes w\mapsto (v,w)\sum_{i=1}^N f_i\otimes f_i,
\end{equation}
where $(f_i)$ is an orthonormal basis of $V$.
It is easy to check that $\Phi(e)^2=N\Phi(e)$. We can now formulate
Brauer's result as follows:

\begin{theorem}\label{brauertheoremo} There exists a surjective
homomorphism $\Phi$ of the algebra $\C_n(N)$, defined over $\CC$
onto $\End_{O(N)}(V^{\otimes n})$. It can be explicitly
described by mapping permutation graphs to the corresponding
permutation of factors in $V^{\otimes n}$, and by mapping
$e_1$ to $\Phi(e)\otimes 1_{n-2}$, where $\Phi(e)$ is as in \ref{Emap}
and $1_{n-2}$ is the identity map of $V^{\otimes n-2}$. Moreover,
this representation is faithful if $N>n$.
\end{theorem}

\subsection{Symplectic case} The same strategy also works
for describing $\End_{Sp(N)}(V^{\otimes n})$, where $V=\CC^N$
is the vector representation of the symplectic group $Sp(N)$, with
$N$ even. Let now $(\ ,\ )$ be a nondegenerate symplectic bilinear
form on $V$, and let $(f_i)$, $(g_i)$ be dual bases with respect to
$(\ ,\ )$. Then we define $\Phi(e)\in\End(V^{\otimes 2})$ by
$\Phi(e)(v\otimes w)=(v,w)\sum_i f_i\otimes g_i$.

\begin{theorem}\label{brauertheoremsp}
There exists a surjective homomorphism $\Phi$ of the algebra $\C_n(-N)$,
defined over $\CC$ onto $\End_{Sp(N)}(V^{\otimes n})$. It can be explicitly
described by mapping permutation graphs to the corresponding signed
permutation of factors in $V^{\otimes n}$, and by mapping
$e_1$ to $-\Phi(e)\otimes 1_{n-2}$, where $\Phi(e)$ is as defined
in this subsection.  Moreover,
this representation is faithful if $|N|>n$.
\end{theorem}

\begin{remark} There also exists a symplectic Brauer algebra which maps
surjectively onto the algebra $\End_{Sp(N)}(V^{\otimes n})$,
and for which the
homomorphism in the last theorem is a little more natural.
We will only use the orthogonal Brauer algebra, although this leads
to some minor complications in notation later.
\end{remark}

\subsection{$q$-Version of Brauer algebra}\label{qversion}
We also need the complex  algebra $\Cc_n(r,q)$ depending on complex
parameters $r$ and $q$. We denote by $\Cc_n$
the corresponding algebra, defined over a field of
rational functions where $r$ and $q$ are considered variables.
The algebra $\Cc_n(r,q)$ is defined via
generators $T_1,\ T_2 \ ...\ T_{n-1}$,
which satisfy the braid relations $T_iT_{i+1}T_i=T_{i+1}T_iT_{i+1}$
and $T_iT_j=T_jT_i$ if $|i-j|>1$ as well as the relations
\vskip .2cm

$(R1)\quad E_iT_i=r^{- 1}E_i$,

$(R2)\quad E_iT_{i-1}^{\pm 1}E_i=r^{\pm 1}E_i,$
\vskip .2cm
where $E_i$ is defined by the equation
$$ (q-q^{-1})(1-E_i)=T_i-T_i^{-1}.\leqno{\indent (D)}$$
\vskip .1cm
The algebras $\Cc_n(r,q)$ have a trace functional $tr$. To define
it, it will be convenient to use the notation
$x=\frac{r-r^{-1}}{q-q^{-1}}+1$. Then $tr$ can
be defined inductively by $tr(1)=1$, $tr(T_i)=r/x$, $tr(E_i)=1/x$
and $tr(a\chi b)=tr(ab)tr(\chi)$
for $a,b\in\Cc_n$ and $\chi\in\{ 1,T_n,E_n\}$.
The algebra $\Cc_n(q^{N-1},q)$ plays the same role for the
Drinfeld-Jimbo quantum group $U_qso_N$ as the Brauer algebra
$\C_n(N)$ plays for the orthogonal group $O(N)$. Moreover,
analogues of  the graph basis for Brauer algebras have been defined
in \cite{GH} (see also references there). In this set-up, one obtains
the algebra $\C_n(N)$ as a limit $\lim_{q\to 1}\Cc_n(q^{N-1},q)$.

\subsection{Algebraic structure}\label{algstruc}
As usual, a Young diagram $\la=(\la_i)$
is an array of boxes, with $\la_i$ boxes in the $i$-th row.
We will freely identify $\la$ with a vector in $\ZZ^m$ whose $i$-th
coordinate is equal to $\la_i$ whenever
$m$ is greater than the number of rows. Let $|\la |$ be the number
of boxes of $\la$.
We denote by $\la'$ the
Young diagram with rows and columns interchanged. In particular, $\la_i'$
is the number of boxes in the $i$-th column.
\vskip .2cm
{\it Structure of $\C_n$ and $\Cc_n$} (a) The algebra $\C_n$  is a direct sum
of full matrix algebras.
Its simple components
are labelled by the Young diagrams
with $n$, $n-2$, $n-4$, ..., 1 resp. 0
boxes.

(b) Let  $U_{n,\la}$ be  a simple $\C_{n,\la}$ module.
The decomposition of  $U_{n,\la}$
into simple $\C_{n-1}$ modules is given by
\begin{equation}\label{restriction}
U_{n,\la}\cong \bigoplus_{\mu} U_{n-1,\mu},
\end{equation}
where the summation goes over all Young diagrams $\mu$ which can
be obtained by either taking away or, if $\la$ has less than $n$
boxes, by adding a box to $\la$.
The labeling of simple components
is uniquely determined by the restriction rule, except for a possible
choice of replacing $\la$ by its transposed $\la'$
simultaneously for all Young diagrams
(see e.g. \cite{WHe}, Lemma 2.11).
The dimension of $U_{n,\la}$ is determined by this formula inductively.
There also exist explicit formulas for them, see e.g. \cite{R}.

(c) The analogous statements for semisimplicity and
restriction rules also  hold for the algebra $\Cc_n$. Here we have the
convention that
the eigenprojection of $T_1$ corresponding to its eigenvalue
$q$ is labeled by the Young diagram $[2]$.

(d) There exists a similar theory
for the integral versions of $\C_n(N)$ and $\Cc_n(r,q)$
defined over suitable rings (see e.g. \cite{GH} and references therein).
It was shown that both series of algebras are cellular
in the sense of \cite{GL}, see \cite{GL} and \cite{Xi}. This means,
in particular, that we can also find suitable bases for the modules
$U_{n,\la}$ with respect to which the generators act via matrices
with entries in those rings. Every simple module of $\C_n(N)$ and $\Cc_n(r,q)$
appears as a quotient of a $U_{n,\la}$ with respect to the annihilator
of a certain bilinear form, also if the algebras themselves are not
semisimple; this can happen for $\C_n(N)$ if $N$ is an integer
or for $\Cc_n(r,q)$ if $q$ is a root of unity and/or if $r=\pm q^k$
for some integer $k$, see \cite{w1}, \cite{w2}.

(e) The labeling of simple components of $\C_n$ and $\C_n(N)$
via Young diagrams
is chosen such that the idempotent $p_{[1^2]}$ corresponds to
the antisymmetrization of $V^{\otimes 2}$ for $O(N)$ and to
the symmetrization of  $V^{\otimes 2}$ for
$Sp(|N|)$ under the map $\Phi$ in Theorems \ref{brauertheoremo} and
 \ref{brauertheoremsp}. This entails that a minimal idempotent $p_\la$
with $\Phi(p_\la)\neq 0$
projects onto a simple $O(N)$-module labeled by $\la$,
and onto a simple $Sp(N)$-module with highest weight $\la'$.

\subsection{Special Quotients}\label{specquot} Let us first do the
orthogonal case.
Let $\Cc_n(M,\ell)$ be the $\CC$ algebra
defined by the $T_i$'s and $E_i$'s as above with $q=e^{\pi i/\ell}$ and
$r=q^{M-1}$. The following statements would hold as well if $q$ is replaced
by any primitive $2\ell$-th root of unity.
We define $\bCcn(M,\ell)$ to be the quotient of  $\Cc_n(M,\ell)$ modulo
the annihilator ideal $I_n(M,\ell)$ of $tr$, i.e.  $I_n(M,\ell)$
consists of all the elements $a\in \Cc_n(M,\ell)$ such that
$tr(ab)=0$ for all $b\in \Cc_n(M,\ell)$.
Then we have the following properties:

(a) We have $I_n(M,\ell)= I_{n+1}(M,\ell)\cap \Cc_n(M,\ell)$, so we
obtain well-defined embeddings $\bCcn(M,\ell)\subset \bar \C_{n+1}(M,\ell)$.

(b) The quotients  $\bCcn(M,\ell)$ are semisimple
and its simple components are
labeled by Young diagrams $\la$ with $n, n-2,\ ...\ $ boxes satisfying
$\la'_1+\la'_2\leq M$ and $\la_1+\la_2\leq \ell+1-M$ for $M>0$.

(c) The restriction rule for simple $\bCcn(M,\ell)$ modules is the same as
for the generic case, as described in the previous chapter, except that
now only the diagrams satisfying the conditions in part (b) are allowed.

For the symplectic case, let $M<0$ be even. We define $\Cc_n(M,\ell)$
and $\bCcn(M,\ell)$ as before with $q=e^{\pi i/2\ell}$ and $r=q^{M-1}$.
Then statements (a) and (c) above also hold in this case, while
for statement (b)  now
the Young diagrams have to satisfy the conditions
$\la_1\leq |M|/2$ and $\la_1'\leq \ell - |M|/2$.

\begin{remark}\label{dimmult} It follows from the semisimplicity
of representations that the multiplicity of a simple $O(N)$ module
$V_\la$ labeled by the Young diagram $\la$ is equal to the dimension
of a simple $\Phi(\C_{n,\la})$-module $\tilde W_{n,\la}$, where
$\Phi$ is the same as in Eq \ref{Emap}. These dimensions can be
calculated inductively by the restriction rule as in \ref{restriction},
but where now only admissible Young diagrams are allowed;
admissible $O(M)$ and $Sp(M)$ diagrams are as defined above
with $\ell=\infty$, i.e. one of the two conditions above becomes void.
Observe that this is equivalent to the tensor product rule for
$V_\la\otimes V$; it is well-known that this is a direct sum of
simple modules $V_\mu$, where $\mu$ runs through all admissible
diagrams which can be obtained by adding or removing a box to/from
$\la$.

While the ring $Gr(G)_\ell$ is not the representation ring of a group,
we will see in the next chapter that it is the representation ring of
a semisimple tensor category $\F=\F(O(N)_\ell)$.
If $X_\la$ is the simple object in $\F$ labeled by $\la\in Gr(G)_\ell$,
and, for  $G=O(N)$ or $G=Sp(N)$,  $X=X_{[1]}$ is the object corresponding to
the vector representation, we again get
that $X_\la\otimes X\cong \bigoplus_\mu X_\mu$, where  $\mu$ runs
through all admissible
diagrams (for $Gr(G)_\ell$)
which can be obtained by adding or removing a box to/from
$\la$. It follows as before that $\End_\F(X^{\otimes n})\cong \bCcn(N,\ell)$.
\end{remark}

\begin{lemma}\label{levelrank1} Let $T_i(r,q)$ and $E_i(r,q)$
denote the generators of $\Cc_n(r,q)$.

(a) The map $T_i(-r^{-1},q)\mapsto -T_i^{-1}(r,q)$
defines an isomorphism between $\Cc_n(-r^{-1},q)$ and $\Cc_n(r,q)$
which maps the simple component $\Cc_n(r,q)_\la$ to
$\Cc_n(-r^{-1},q)_{\la'}$.

(b) Let $q$ be a primitive $2\ell$-th root of
unity, and let $1<M<\ell$. Then the map in (a) also induces
an isomorphism between $\bCcn(M,\ell)$ and
$\bCcn(\ell+2-M,\ell)$. Moreover,
these isomorphisms are compatible with the standard inclusions
$\bCcn(M,\ell)\subset \bar\C_{n+1}(M,\ell)$.

(c) Let $q$ be a primitive $4\ell$-th root of unity, and let $-2\ell<M<-1$,
with $M$ even.
Then the map in (a) also induces an isomorphism between
$\bCcn(M,\ell)$ and $\bCcn(-2\ell+|M|+2,\ell)$, which is again compatible
with inclusions.
\end{lemma}

$Proof.$ The isomorphism property is easy to check.
It maps the eigenprojection of $T_1(-r^{-1},q)$ for $q$ (the component
of $\Cc_2(-r^{-1},q)$ labeled by $[2]$) to the eigenprojection of
$T_i(r,q)$ of $-q^{-1}$, which is labeled by the Young diagram $[1^2]$.
By uniqueness of labeling of the simple components, this implies
that the simple component $\Cc_n(-r^{-1},q)_\la$ is mapped to
the simple component $\Cc_n(r,q)_{\la'}$. This shows (a).

Now observe that the map on $\Cc_n(-r^{-1},q)$ induced by the concatenation
of our map with
the Markov trace on $(\Cc(r,q))$ coincides with the Markov
trace on  $\Cc_n(-r^{-1},q)$. Indeed, by the inductive definition of
the trace, it suffices to check this  for the generators. As $x(r,q)=
x(-r^{-1},q)$ (see definitions of $(\Cc(r,q))$), it follows that
$tr(E_i(r,q))=tr(E_i(-r^{-1},q))$. Using relation $(D)$, one also
sees that $tr(T_i^{-1}(r,q))=r^{-1}/x$. Hence also $tr(T_i(-r^{-1}))=
tr(-T_i^{-1}(r,q))$, which proves our claim.

Now if $q$ is a primitive
$2\ell$-th root of unity, and $r=q^{M-1}$, we have $-r^{-1}=q^\ell q^{1-M}
= q^{\ell+1-M}$. It follows from the two previous paragraphs that we
have a trace preserving isomorphism between $\Cc_n(M,\ell)$ and
$\Cc_n(\ell+2-M,\ell)$. The claim follows from this and the definitions
of $\bCcn(M,\ell)$ and  $\bCcn(\ell+2-M,\ell)$. The proof for (c)
goes similarly.

\section{Grothendieck semirings}\label{Grothendiecksemi}

The main purpose of this section is to prove Prop. \ref{fusionsym},
which is a case of what is known as level-rank duality. This has
been observed before in the physics and mathematics literature
(see e.g. the discussions to Chapters 16 and 17 in
\cite{FMS} and e.g. the papers \cite{BB} and \cite{MPS}).
We will give a fairly simple and elementary proof
for the particular cases which we need, which is more or less the
same approach as the one in the papers \cite{BB} and \cite{MPS},
which go back to results in \cite{w1} and \cite{TW2}.
It will hopefully make the translation between the various
approaches to our setting as well as fixing notations easier.
It should also explain
the rings $Gr(G)_\ell$ of the first section to some extent,
and will put some of the other topics discussed later into
a more conceptual setting.

\subsection{Semirings from idempotents}
Let $G=O(N)$ or $G=Sp(N)$, and let $\Phi : \C_n(\pm N)\to
\End_G(V^{\otimes n})$ be the surjective homomorphism
mentioned in Section \ref{defsec}.
As an idempotent $p\in \C_n(N)$
corresponds to the subrepresentation
$\Phi(p)V^{\otimes n}$, we can translate the tensor product structure of
$Rep(G)$ into the setting of the algebras $\C_n(N)$. More generally,
this can also be done as well for the algebras $\C_n$ and $\Cc_n$ which will
lead to the definition of a formal Grothendieck semiring $\Grin_+$.
The following definitions and constructions have
appeared before, e.g. in $K$-theory for $C^*$-algebras and
the idempotent construction for categories, also sometimes referred to
as the Karoubian.

We say that two idempotents $p$ and $q$ in an algebra $\A$ are
{\it conjugation equivalent}
if there exist elements
$u$ and $v$ in $\A$ such that $p=uv$ and $q=vu$. More generally, we say that
two idempotents $p\in\C_n$ and $q\in\C_m$ are equivalent if we can find
nonnegative integers $n_1$ and $n_2$ such that
$(\frac{1}{ x}e)^{\otimes n_1}\otimes p$ and
$(\frac{1}{ x}e)^{\otimes n_2}\otimes q$ are conjugation equivalent
in the algebra $\C_{n+2n_1}=\C_{m+2n_2}$. For a given idempotent $p\in\C_n$
we denote by $[p]$ its equivalence class. It is easy to check
that for idempotents $p\in\C_n$ and $q\in\C_m$ we obtain an idempotent
$p\otimes q\in\C_{n+m}$ which gives rise to a well-defined multiplication
$[p][q]=[p\otimes q]$.

Addition is defined as follows: Assume $p$ and $q$ are idempotents
in $\C_n$ such that $q$ is conjugation equivalent to an idempotent
$q'\in\C_n$ satisfying $pq'=0=q'p$. Then also $p+q'$ is an idempotent,
and we define $[p]+[q]=[p+q']$. It is easy to check that this is
well-defined. Moreover, it is known that for fixed
Young diagram $\la$ the dimension of the simple components
$\C_{|\la|+2k,\la}$ goes to infinity if $k\to\infty$ (see the
formula in the last section,(b)).
Hence, by tensoring $p$ and $q$ by a projection $(\frac{1}{ x}e)^{\otimes k}$
for suitably large $k$, if necessary, we can always assume the existence
of a $q'$ as above; indeed, it is easy to see that such a $q'$ exists
whenever the sum of the ranks of $p$ and $q$ in an irreducible
representation is less than its dimension.

Let $p_\la$ and $p_\mu$ be minimal idempotents
in the simple components of $\C_n$ and $\C_m$
labeled by $\la$ and $\mu$ respectively. We shall often just
write $[\la]$ and $[\mu]$ for $[p_\la]$ and $[p_\mu]$.
If we define $d^\nu_{\la\mu}$ to be the rank of the idempotent
 $p_\la\otimes p_\mu$ in the irreducible representation of
$\C_{n+m}$ labeled by the Young diagram $\nu$, we obtain
\begin{equation}\label{multiplication}
[\la][\mu]=[p_\la\otimes p_\mu]=\sum_\nu d_{\la\mu}^\nu [\nu].
\end{equation}
We denote the semiring of equivalence classes of idempotents
in $\cup \C_n$ with addition and multiplication as defined above by $\Grin$.
Observe that $\Grin$ can also be defined as the free $\N$ module
spanned by
the equivalence classes of minimal idempotents in  $\cup \C_n$,
with multiplication defined by Eq. \ref{multiplication}.

Both definitions work as well
for the sequence of algebras $\C_n(N)$ and for the sequence
of algebras $\Cc_n$. In the latter case
the tensor product $a\otimes b$ of elements
$a\in\Cc_n$ and $b\in\Cc_n$ is defined by
 $a\otimes b = a \sh_n(b)$, where $\sh_n :\Cc_m\to\Cc_{n+m}$
is the homomorphism defined by $\sh_n(T_i)=T_{i+n}$,  $\sh_n(E_i)=E_{i+n}$,
for $i=1,2,\ ...\ m-1$.

It should be clear that the same constructions also work for
the sequences of algebras $\Cc_n(r,q)$ and $\bCcn(n,\ell)$.
Similarly, we can also define a semiring
for the sequence of
group algebras $\CC S_n$ of the symmetric groups $S_n$, $n\in\N$.
In view of Schur duality, a minimal idempotent $p_\la\in \CC S_n$
 corresponds to an irreducible representation $F^\la$
of $Gl(N)$ for $N$ sufficiently large.
The resulting semiring will be denoted by $Gr(Gl(\infty))$.
We denote the structure coefficients in the semiring
$Gr(Gl(\infty))$ by $c_{\la\mu}^\nu$;
they give the multiplicity of the simple $Gl(N)$-module
$F^\nu$ in the tensor product $F^\la\otimes F^\mu$.

\begin{definition}\label{semirings}
We define  $Gr(\C(N))_+$ resp $Gr(\Cc(r,q))_+$ resp $Gr(\bCcn(N,\ell))$
to be the semirings derived from the  equivalence
classes of idempotents of the algebras  $\C_n(N)$,
resp $\Cc_n(r,q)$, resp $\bCcn(N,\ell)$
$n\in\N$, $N\in\ZZ$, and $r,q\in\CC$ and we denote by
$Gr(Gl(\infty))_+$ the semiring obtained from the group algebras
$\CC S_n$ of the symmetric groups $S_n$, $n\in\N$.
The corresponding rings are denoted by the same symbols, without
the $+$.
\end{definition}

\begin{remark} It is possible to construct tensor categories from
the algebras mentioned above, see \cite{TW2}, \cite{BB}, \cite{TbW}.
The just defined semirings would be the Grothendieck semirings
of these tensor categories. We will discuss this and other constructions
of such categories below.
\end{remark}

\subsection{Algebraic properties}

\begin{lemma}\label{genericring} Let $d_{\la\mu}^\nu(N)$ and
$d_{\la\mu}^\nu$ be the structure coefficients for $Gr(SO(N))_+$ and
$\Grin$, see Section \ref{Grin}. Then we have

(a)  $d_{\la\mu}^\nu(N)=d_{\la\mu}^\nu$ for $|N|$ sufficiently large,

(b) $Gr(\C)_+\cong \Grin_+$
\end{lemma}

$Proof.$ Let $p_\la, p_\mu$ be minimal idempotents in $\C_n(N)$,
defined over $\Q(x)$. If $N$ is not one of the finitely many
poles in the coefficients of $p_\la$, $p_\mu$ (with respect to the
graph basis), we also obtain well-defined idempotents $p_\la(N)$,
$p_\mu(N)$ in the $\Q$-algebra $\C_n(N)$. Obviously,
both $p_\la\otimes p_\mu$ and  $p_\la(N)\otimes p_\mu(N)$
act as idempotents with the same rank $d^\nu_{\la\mu}$
on the module $U_{n,\nu}$. On the other hand, for $N>n$,
$\End_{O(N)}(V^{\otimes n})\cong \C_n(N)$, and hence $d^\nu_{\la\mu}$
is equal to the multiplicity  $d_{\la\mu}^\nu(N)$
of the simple $O(N)$-module $V_\nu$ in
$V_\la\otimes V_\mu$.
Finally, assuming $N$ big enough, we can identify the Young diagram
$\la$ as well with the corresponding $SO(N)$ dominant weight.
This proves both (a) and (b).

\begin{lemma}\label{transp}
Let $\la'$ denote the transpose of the Young diagram $\la$.

(a) The map
$[\la]\mapsto [\la']$ defines an automorphism of $\Grin$.

(b)  The map
$[\la]\mapsto [\la']$ defines an isomorphism between  $Gr(\bCcn(N,\ell))$
and  $Gr(\bCcn(\ell + 2-N,\ell))$
\end{lemma}

$Proof.$  This follows from Lemma \ref{levelrank1}.
$\Box$

\subsection{Ribbon tensor categories}
We give here some very limited background on such tensor categories.
The reader can find very good introductions to tensor categories
suitable for our context in \cite{Kassel} and \cite{Turaev}.

In the following we assume $\T$ to be a  monoidal rigid ribbon
tensor category. This means that we have a functor $\otimes : \T\times \T
\to \T$ satisfying certain associativity conditions (monoidal) and
rigidity, i.e. every object has a dual object (in the sense of a dual
representation for group representations). We also have a braiding
which means that there exist for any pair of objects $X,Y$ in $\T$
a canonical morphism $c_{X,Y}: X\otimes Y \to Y\otimes X$ subject to
several axioms. In particular, we obtain for any object $X$ a
representation of Artin's braid group $B_n$ into $\End(X^{\otimes n})$.
Conversely, one can construct from suitable series of braid group
representations a monoidal rigid ribbon tensor categories.
This has been carried out for the algebras $\bCcn(N,\ell)$ in \cite{TW2}.
In particular, $Gr( \bCc(N,\ell))_+$ is the Grothendieck semiring
of this category.

A special class of ribbon categories are modular tensor categories.
To explain  it, we first remark that one can define for any object
$X$ in a rigid ribbon category a trace functional $Tr_X$ on
$\End(X)$ (see e.g. \cite{Kassel} Ch. XIV, 4). In particular,
one defines the categorical dimension $d(X)$ of an object $X$
by $d(X)=Tr_X(1_X)$, where $1_X$ is the identity in $\End(X)$.
If $X_\la,X_\mu$ are simple objects in our tensor category,
we define the scalar
$s_{\la\mu}=Tr_{X_\la\otimes X_\mu}
(c_{X_\mu\otimes X_\la}\circ c_{X_\la\otimes X_\mu})$.
A modular tensor category is defined to be a semisimple rigid ribbon
tensor category with only finitely many simple objects up to
isomorphism for which the matrix $S=(s_{\la\mu})$ is invertible.
Now for any simple object $X_\mu$ of a ribbon category $\T$ with
nonzero dimension $d(X_\mu)$ one can show the following
 Verlinde formula (see e.g. \cite{Turaev}, Th. 4.5.1):
The map $\phi_\mu: [X_\la]\mapsto s_{\la\mu}/d(X_\mu)$
defines a multiplicative character of $Gr(\T)$. In particular,
for a modular tensor category, we thus obtain a complete set of
characters of the abelian ring $Gr(\T)$; indeed, the characters
$(\phi_\mu)$ are linearly independent due to the invertibility
of the $S$-matrix.

\subsection{Examples}
An important class of ribbon tensor
categories is provided by the Drinfeld-Jimbo
$q$-deformations $U_q\g$ of the universal enveloping algebra of a
semisimple Lie algebra $\g$.  For generic $q$ (i.e. $q$ not a root of unity
except for $q=\pm 1$) the category of finite dimensional representations
of $U_q\g$ is semisimple. Moreover, its simple representations are
labeled by the dominant weights of $\g$. This provides an
isomorphism between the Grothendieck semirings of $Rep(\g)$ and of
$Rep(U_q\g)$, where $Rep$ refers to finite dimensional representations.

It is more complicated to construct modular tensor categories.
We will rely here on constructions using tilting modules of quantum
groups at roots of unity (see \cite{A}, \cite{AP}). There are also
more elementary constructions (see \cite{TW2}, \cite{Bl},\cite{BB}) which
would at least cover some of the cases needed here; but they
are not as general. Additional constructions
can be found in \cite{KL} and \cite{Wa}.
We have the following theorem:

\begin{theorem}\label{AndersenP} For each connected
and simply connected Lie group
$G$ and positive integer $\ell$  exists a modular tensor category
$\F=\F(G_\ell)$
whose Grothendieck semiring is isomorphic to $Gr(G)_{\ell, +}$ as in
Section \ref{generalconstruction}.
\end{theorem}

$Proof.$ This is essentially contained in \cite{A}, \cite{AP}.
The authors construct a braided tensor category $\F$
as a quotient of the category of tilting modules of a quantum group
$U_q\g$ for $q$ a root of unity. Its simple objects $T_\la$
(up to isomorphism) are labeled by the same index set as
the basis elements of $Gr(G)_\ell$, where $G$ is the
simply connected Lie group corresponding to $\g$ and $\ell$ is
determined by the degree of the root of unity. In particular, the entries of
the $S$-matrix can be calculated using the theory of $R$-matrices
(see e.g. \cite{TW1}, Lemma 3.5.1).
If $q=e^{\pi i/\ell}$, one can explicitly calculate that
$s_{\la\mu}=\chi^\la((\mu +\rho)/\ell)d(X_\mu),$
where $\chi^\la$ is given by the Weyl character formula,
and $d(X_\mu)=d(\mu)$  as in Section \ref{generalconstruction};
in particular $(s_{\la\mu})$ is invertible.
Hence the characters $\phi_\mu$ mentioned in this section
coincide with the ones in Proposition \ref{fusionring},
and the Grothendieck semiring of $\F$ is isomorphic to $Gr(G)_\ell$.
In particular, the subset $Gr(G)_{\ell, +}\subset Gr(G)_\ell$
(see Remark \ref{Grplus}) is closed under multiplication and hence
forms a semiring.

\subsection{Identification with  $Gr(\bCcn(N,\ell))$} As outlined
in the beginning of this section, we can identify $Gr(O(N))$ with
$Gr(\Phi(\C(N)))$, as $\Phi(\C_n(N))\cong \End_{O(N)}(V^{\otimes n})$.
We are going to prove similarly that $Gr(\bCcn(N,\ell))$ is isomorphic
to the Grothendieck semiring of a suitable fusion category. The only
slight subtlety comes from the fact that  $Gr(\bCcn(N,\ell))$ would
correspond to a fusion category in connection with $O(N)$, while the
ones obtained from quantum groups would correspond to $SO(N)$.

\begin{prop}\label{identification}
(a) $Gr(\bCc(-N,2\ell))_+\cong Gr(Sp(N))_{\ell,+}$ for $N$ even.

(b) $Gr(\bCc(N,\ell))_+\cong Gr(O(N))_{\ell,+}$
\end{prop}

$Proof.$ We will check that the stated isomorphisms are given by
associating the equivalence class $[p_\la ]$, with $p_\la$
a minimal idempotent in $\Cc_{n,\la}(q^{\pm N-1},q)$, to the
element $\la$ in $Gr(O(N))_{\ell,+}$ resp. $Gr(Sp(N))_{\ell,+}$.
Let $\F$ be one of the fusion categories related to an
orthogonal or symplectic group, and
let $X=X_{[1]}$ be a simple object
in $\F$ whose isomorphism class is labeled by the Young diagram $[1]$.
Then it is well-known that the representation of the braid group
$B_n$ derived from the braiding morphisms $c_{X,X}$
is a quotient of $\Cc(q^{N-1},q)$ (for $G=O(N)$)
or of $\Cc(q^{-1-N},q)$ (for $G=Sp(N)$), see e.g. \cite{Tur}, \cite{w1}
or \cite{TbW}. Moreover, the pull-back of the categorical trace
$Tr_{X^{\otimes n}}$ on $\End(X^{\otimes n})$ induces the Markov
trace of $\Cc_n$, up to a multiple (see e.g. \cite{TbW}, Cor 5.3).
Hence the quotient
$\bCcn(q^{N-1},q)=\Cc_n(q^{N-1},q)/Ann(tr)$ is also a quotient
in the image of $B_n$ in  $\End(X^{\otimes n})$. Moreover,
if we choose for the required $2\ell$-th (orthogonal case)
or $4\ell$-th (symplectic case) root of unity
$q=e^{\pi i/\ell}$ or $q=e^{\pi i/2\ell}$, respectively,
$Tr_{X^{\otimes n}}$
induces an inner product on  $\End(X^{\otimes n})$
(see \cite{W5}). Hence the image of $Ann(tr)$ in $\End(X^{\otimes n})$
is 0, and the aforementioned quotient is a subalgebra
in  $\End(X^{\otimes n})$.

Comparing the tensor product rules in $\F$ with the embeddings of
$\bCcn\subset \bCc_{n+1}$, one checks by induction on $n$
that for $G=Sp(N)$ and for $G=O(N)$ with $N$
odd, we have $\End_\F(X^{\otimes n})\cong \bCcn(q^{-N-1},q)$ (for $Sp(N)$),
and $\cong \bCcn(q^{N-1},q)$ (for $O(N)$), see Sections
\ref{algstruc}, \ref{specquot} and Remark \ref{dimmult}. Moreover,
these isomorphisms are compatible with the tensor product operation.
Hence we get the desired isomorphism of Grothendieck semirings.

For the even-dimensional case, we will show that $\bCcn(q^{N-1},q)$
is isomorphic to a subalgebra of $\End(X^{\otimes n})$ for which
the inclusion is compatible with the restriction rules from
$Gr(O(N))_\ell$ to $Gr(SO(N))_\ell$. It follows from Brauer's results
that $\Phi(\tilde p_\la)V^{\otimes n}$ decomposes as an $SO(N)$-module
according to the restriction rules, as e.g. stated in Section
\ref{fullortho}; here $\tilde p_\la$ is a minimal idempotent in $\C_{n,\la}(N)$.
We claim that it suffices to find, for given $\la\in D_\ell$, an idempotent
$p_\la=p_\la(q)\in \Cc_{n,\la}(q^{N-1},q)$ such that $p_\la(1)=\tilde p_\la$
and which is well-defined at $q=e^{\pi i/\ell}$.
To see this, consider the explicit representation of
$\Cc_{n,\la}(q^{N-1},q)$ into $\End(V^{\otimes n})$ derived from
Jimbo's  $R$-matrices, see \cite{w1}, Sect. 5, or \cite{Tur}.
Here the basis vectors of $V=\C^N$ are weight vectors of the quantum
group $U_q{\mathfrak so}_N$ for any $q$. As the image of
$\Cc_{n,\la}(q^{N-1},q)$ commutes with the action of  $U_q{\mathfrak so}_N$,
each weight space of $V^{\otimes n}$ is invariant under
$\Cc_{n,\la}(q^{N-1},q)$ for any $q$. In particular, the
$U_q{\mathfrak so}_N$ character of
$p_\la V^{\otimes n}$ does not depend on $q$, and the restriction rule
also holds for $p_\la(e^{\pi i/\ell})$. As $Gr(O(N))_\ell$ was defined
from $Gr(SO(N))_\ell$ via these restriction rules, the claim would
follow.

Inductive formulas for certain minimal idempotents,
called path idempotents, were defined in \cite{RW}, using
prior  work in \cite{WHe} for Hecke algebras;
see also \cite{lr} for a similar approach. The formulas are written
down explicitly, involving quantities such as hook lengths in
Young diagrams. One can check that they are well-defined
for paths only involving diagrams in the labeling set of $Gr(O(N))_\ell$
for $q=e^{\pi i/\ell}$, and that one can get a partition of unity in
$\bCc_n(q^{N-1},q)$ in terms of such idempotents.
This is fairly straightforward, using the explicit formulas in
\cite{RW}. This finishes the proof.

\medskip
{\it Proof of Prop. \ref{fusionsym}}: In the orthogonal case,
it follows from Prop. \ref{identification} that $Gr(O(N))_{\ell,+}\cong
Gr\bCc(N,\ell)_+$. Hence $Gr(O(N))_{\ell,+}\cong Gr(O(\ell + 2-N))_{\ell,+}$,
by Lemma \ref{transp}. The proof goes the same in the symplectic case.

\section{Restriction Coefficients}

Let $p_\la$ be a minimal idempotent in $\CC S_n$. Then it can be written
as a  sum of mutually commuting minimal idempotents in $\C_n$.
Let $b^\la_\mu$ be
the number of those idempotents which are in $\C_{n,\mu}$; equivalently,
$b^\la_\mu$ is the trace of $p_\la$ in an irreducible $\C_{n,\mu}$-module.
Again, these coefficients depend on whether we take the orthogonal or
symplectic Brauer algebra, which will be made precise below.

The connection of these coefficients with representation theory of Lie groups
is easily established as follows.
Let $V$ be an $N$-dimensional vector space with $N>n$. By Schur duality,
$F^\la=p_\la V^{\otimes n}$ is an irreducible $Gl(N)$-module. As $N>n$,
we obtain a faithful representation of the Brauer algebra $\C_n(N)$ on
$V^{\otimes n}$ such that its image is isomorphic to $\End_G(V^{\otimes n})$
for $G=O(N)$ or $G=Sp(N)$. Hence it follows that the coefficient
$b^\la_\mu$ determines the multiplicity of the simple $G$-module $V_\mu$
in $F^\la$. If $N\leq n$, this multiplicity will be denoted by $b^\la_\mu(G)$.
Formulas for the coefficients $b^\la_\mu$ and  $b^\la_\mu(G)$
already appeared in \cite{Wy} (see Theorems 7.8F and 7.9C).
An explicit formula for  $b^\la_\mu$ was computed by Littlewood
\cite{Li} as follows,
where part (c) is a simple consequence of parts (a) and (b) as well as
Lemma \ref{transp}. See also \cite{KT}, Section 1.5
for a more detailed description
of this result and its proof, and \cite{EW} whose general proof
also includes a derivation of the classical case.

\begin{theorem}\label{littlew}
(Littlewood) In the orthogonal case,
we have $b_\mu^\la=\sum_\nu c^\la_{\mu\nu}$,
where the summation goes over all Young diagrams which have an even number
of boxes in each row, and $c^\la_{\mu\nu}$ is the multiplicity of the
simple $Gl(N)$-module $F^\nu$ in $F^\la\otimes F^\mu$.
In particular, $b^\la_\mu=0$ unless $\mu\subset\la$, and $b^\la_\la=1$.

(b) In the symplectic case, we have $b_\mu^\la=\sum_\nu c^\la_{\mu\nu}$,
where the summation goes over all Young diagrams which have an even number
of boxes in each column, and $c^\la_{\mu\nu}$ is as in (a).

(c) The orthogonal and symplectic restriction coefficients are related by
$b^\la_\mu(O)=b^{\la'}_{\mu'}(Sp)$.
\end{theorem}

Recall that for  $G=O(N)$ or $Sp(N)$ we defined a labeling set
$D(G)$ as a certain subset of the set of Young
diagrams (see Sections \ref{sympsect} and
\ref{fullortho}). Recall that we have defined an action of
a reflection group $W$ on Young diagrams for which the closure of $D(G)$
is a fundamental domain, see Sections \ref{lrefl} and \ref{sympsect}.
The corresponding sets $D(G)'=\{\la', \la\in D(G)\}$ are
the usual labeling sets of irreducible representations of $G$,
see Theorem \ref{classicsym}.
Let $b^\la_\mu$ be the multiplicity of the simple $G$-module
$V_\mu$ in the simple $Gl(N)$-module $F^\la$.
We can now express the restriction coefficients $b^\la_\mu(G)$ as follows
(see also \cite{EW}, \cite{KT}, \cite{Ki}, \cite{Su}).

\begin{theorem}\label{result} With notations as above,
the restriction multiplicity $b^\la_\mu(G)$  is given by
$$b^\la_\mu(G)=\sum_{w\in W} \ve(w)b^\la_{w'.\mu}.$$
\end{theorem}

$Proof.$
Let  $p_\la\in\CC S_n$ be a minimal idempotent, and let
$\I(G)$ be the ideal in $\Grin$ such that $Gr(G)\cong \Grin/I(G)$.
Then we have
$$[p_\la] = \sum_\mu b^\la_\mu [\mu] \equiv
\sum_{\mu\in D(G)'} (\sum_w \ve(w)b^\la_{w'.\mu})\ [\mu ]\quad
{\rm mod} \I(G),$$
by Theorem \ref{classicsym}.
$\Box$

\begin{remark}  It is tempting to consider
the question whether there exist analogs of restriction coefficients
for fusion categories. This does not seem to be possible with
the categories we have considered here, e.g. there is no
inclusion of quantum groups $U_q{\mathfrak so}_N\subset U_q{\mathfrak sl}_N$.
However, it is expected that one can extend fusion categories of type $A$
by additional objects corresponding to representations of
a twisted Kac-Moody algebra. Tensoring a twisted
with an untwisted representation would formally correspond to
a restriction-induction process. In one of the examples, the
multiplicities would be modifications of restriction coefficients from
$Gl(N)$ to $O(N)$, which, for given diagrams and
large enough $\ell$ would be equal
to the classical coefficients. We plan to study these questions
in connection with subfactors of von Neumann factors.
\end{remark}

\section{Structure of Brauer algebra and generalizations}

\subsection{Introduction}
When the first version of this paper was written and submitted, it was
still an unsolved problem how to determine the structure of Brauer algebras
in the case when they are not semisimple. Our work so far, in combination
with the work \cite{So} on tilting modules,
suggested a conjecture what the solution
should be, at least as far as determining the dimensions of the simple
modules was concerned.
This conjecture was also obtained, independently,
in \cite{CdVM} and, more recently, a proof has been posted in \cite{Mt}.
Moreover, the  paper \cite{Hu} has appeared in the meantime
from which one can determine, in connection with \cite{So},
the structure of  $\Cc_n(q^{-N-1},q)$ with $N$ even
(the symplectic case) in many cases. We extend this to the generic
case (i.e. for all but finitely many $q$) along the approach
outlined in the earlier version.
This might also give a more conceptual
understanding of the situation in the classic Brauer algebra case.

\subsection{Tilting modules and Kazhdan-Lusztig polynomials}
We will need to review some basic facts of Kazhdan-Lustig polynomials. See
\cite{Hm} and, in particular, \cite{So1} for introductions
to this subject. If $W$ is a reflection group, then Kazhdan-Lusztig
polynomials $P_{v,w}$ are defined for any pair $w,v$ of elements
in $W$ via induction on the length of $w$. In particular, $P_{v,w}\neq 0$
only if $v\leq w$ in Bruhat order. Similarly, if $U\subset W$ is a parabolic
subgroup, one can define parabolic Kazhdan-Lusztig polynomials
$P_{\bar v,\bar w}$ for the left cosets $W/U$.
We will need this for the inclusion $W\subset \Wl$, where $W$ is
a Weyl group and $\Wl$ is an affine reflection group,
see Section \ref{generalconstruction}. In this case, there exists a
well-known geometric
model for the left cosets $\Wl/W$ in terms of alcoves in
the dominant Weyl chamber. In particular, the parabolic Kazhdan-Lusztig
polynomials can be calculated via paths from the fundamental alcove
to the given one, see e.g.
\cite{GW2} for a review and for further references.  We will need
the following definition, with notations as
in Section \ref{generalconstruction}.

\begin{definition}\label{defamu}
Let $\la,\mu$ be dominant weights such that $w.\la\neq \la$ for
any nontrivial $w\in\Wl$  Assume that there
exists $w\in\Wl$ such that $\nu=w^{-1}.\la\in D_\ell$. Then we define
$a^\la_\mu=P_{\bar v,\bar w}(1)$ if $v.\nu=\mu$, and $a^\la_\mu=0$ if $\nu$
is not in the $\Wl$-orbit of $\mu$.
\end{definition}

To distinguish between the various cases,
we will write $a^\la_\mu(M,\ell)$ for $G=SO(M)$ (where we assume
$M$ to be odd) and $a^\la_\mu(-M,\ell)$
for $G=Sp(M)$, with $M$ even.
Observe that we  also have actions of $S_\infty\subset W=W(D_\infty)$,
(or of $W=W(B_\infty)$ in the symplectic
case) on infinite sequences as defined in Section \ref{lrefl}.
We then define similarly
the quantity $a^\la_\mu(\pm N)$ for Young diagrams $\la,\mu$ via
parabolic Kazhdan-Lusztig polynomials for $S_\infty\subset W$
with the action $w.\la$ defined in Section \ref{lrefl};
here $a^\la_\mu(-N)$ is defined for the symplectic case $G=Sp(N)$, see Section
\ref{sympsect}.
We now have the following simple lemma, which can be considered a
refinement of Lemma \ref{redlem}.

\begin{lemma}\label{KLlimit} We have $a^\la_\mu(M,N+M-2)=a^\la_\mu(N)$
and  $a^\la_\mu(-M,(N+M+2)/2)=a^\la_\mu(-N)$ independent of $M$ provided
it is sufficiently large.
\end{lemma}

$Proof.$ We assume that $\la$ and $\mu$ are in the same $W$ orbit;
otherwise there is nothing to show. It follows from the discussion
in Section \ref{lrefl} and, in particular, the algorithm in the proof
of Lemma \ref{redlem} that for large enough $m=M/2$, the subgroup
of $\Wl(m)$ generated by the reflections needed to map $\la$ and $\mu$
into $\DN$ is independent of $m$ and can be viewed as a subgroup
of $W=W(D_\infty)$. Hence also the corresponding
Kazhdan-Lusztig polynomial is independent of $m$, and coincides
with the one calculated within $W$. The same argument works
in the symplectic case.

\medskip

Tilting modules were defined for quantum groups by H.H. Andersen \cite{A}
and before, for algebraic groups, by Donkin.
A tilting module has a filtration in terms of Weyl modules; this implies,
in particular, that its character is a sum of ordinary simple characters.
More precisely, it has been shown by
Soergel in \cite{So} that the character $\chi^\la_T$
of an indecomposable tilting
module $T_\la$ of the quantum group $U_q\g$ with highest weight $\la$
can be written as
\begin{equation}\label{soergeltilt}
\chi_T^\la=\sum_\mu a^\la_\mu \chi^\mu,
\end{equation}
where $\chi^\mu$ is the character of the simple $\g$-module with
highest weight $\mu$ and $a^\la_\mu$ is as in Def. \ref{defamu};
here $\ell$ is related to the degree of the root of unity $q$
in $U_q\g$ as indicated in our examples.
We also note here that the vector representation
$V$ of $\U$, with $\g$ of classical Lie type, is a tilting module,
and so is any tensor power $V^{\otimes n}$ and any direct
summand of $V^{\otimes n}$. Finally, any tilting
module is a direct sum of indecomposable tilting modules, and there
exists exactly one indecomposable tilting module $T_\la$
with highest weight $\la$ for each dominant weight $\la$.

\subsection{$q$-Brauer algebras}
As already stated before,
the algebras $\Cc_n(q^{-N-1},q)$ map surjectively
onto $\End_{\bf U}(V^{\otimes n})$, where
${\bf U}=U_q{\mathfrak  sp}_N$, $N$ even,
for $q$ not a root of unity. Moreover, for $N>n$, this representation
is faithful (see Theorems \ref{brauertheoremo} and \ref {brauertheoremsp},
which carry over to the $q$-deformations).
Recall that
for $Sp(N)$ we associate to the dominant weight $\la$ the Young diagram
$\la'$, which has $\la_i$ boxes in its $i$-th $column$, for
labeling the corresponding $\Cc_n(q^{-N-1},q)$ module.
It is shown in \cite{Hu} that this is also true
at roots of unity in context of tilting modules of $U_q{\mathfrak  sp}_N$.
So if $T_\la\subset V^{\otimes n}$ is an indecomposable tilting module
with highest weight $\la$, we can find a minimal projection
$p_{\la'}\in\Cc_n(q^{-N-1},q)$ such that $T_\la=p_{\la'} V^{\otimes n}$.
Moreover, the equivalence classes of minimal idempotents in
$\Cc_n(q^{-N-1},q)$ are in 1-1 correspondence to
indecomposable tilting modules
(again: any minimal idempotent in $\Cc_n(q^{-N-1},q)$ projects
onto a direct summand of the tilting module $V^{\otimes n}$, which hence
is itself an indecomposable tilting module).
We have the following simple consequence of the theorems of
Soergel and Hu:

\begin{theorem}\label{tiltdual} Assume that $M>n$ and $M$ even. Let
$d_{n,\la}(-M,\ell)$ be the dimension of a simple $\Cc_{n,\la}(q^{-M-1},q)$-
module, with $q$ a primitive $4\ell$-th root of unity,
and let $d_{n,\mu}$ be the dimension of $U_{n,\mu}$, see
Section \ref{algstruc}. Then we have
$\sum_\la d_{n,\la'}(-M,\ell)a^\la_\mu(-M,\ell)=d_{n,\mu'}$. These
equations completely determine the dimensions $d_{n,\la}(-M,\ell)$
of all simple $\Cc_{n,\la}(q^{-M-1},q)$ modules.
\end{theorem}

$Proof.$ Decomposing $V^{\otimes n}$ into a direct sum of
simple $Sp(M)$-modules, and into a direct sum of indecomposable
tilting modules, and comparing characters,
with our labeling conventions, we get
$$\chi_{V^{\otimes n}}\ =\ \sum_\mu d_{n,\mu'}\chi^\mu
= \sum_\la d_{n,\la'}(-M,\ell)\chi^\la_T.$$
The claim now follows from expanding $\chi^\la_T$
in ordinary characters, see Eq. \ref{soergeltilt}, and comparing the
coefficients of $\chi^\mu$ on both sides. Finally, it is well-known
that $a^\la_\mu(-M,\ell)$ is nonzero only if $\mu\leq \la$ in Bruhat order.
Hence we get a triangular system of equations for the dimensions
$d_{n,\la}(-M,\ell)$ in terms of the generic dimensions $d_{n,\la}$,
which can be computed inductively, see Section \ref{algstruc}.

\begin{remark}\label{decompmat}
1. The condition $n<M$ is  necessary to get a faithful representation
of $\Cc_n(q^{-M-1},q)$ in $\End(V^{\otimes n})$ (see \cite{Hu}).

2. The matrix $(a^{\la'}_{\mu'}(-M,\ell))$ is usually called
the decomposition matrix for  $\Cc_{n}(q^{-M-1},q)$,
where $q$ is a primitive $4\ell$-th root of unity.
It is not hard to show that $a^{\la'}_{\mu'}(-M,\ell)$ is the multiplicity
of the simple $\Cc_{n}(q^{-M-1},q)$-module labeled by $\la$
in the module $U_{n,\mu}$. This second interpretation
can be fairly easily proven if one takes for $U_{n,\mu}$
the space of all highest weight vectors  of weight $\mu'$
in $V^{\otimes n}$ for a suitable integral form
of $U_q{\mathfrak sp}_{2m}$, using Lusztig's canonical bases,
see \cite{Hu}.
\end{remark}

\subsection{Decomposition matrices for $\Cc_{n}(q^{-N-1},q)$, $N$ even}
The methods in the last section only work for the image of the
representation of  $\Cc_{n}(q^{-M-1},q)$ in $\End(V^{\otimes n})$.
Using level-rank duality, Lemma \ref{levelrank1}, we can extend this
to the full algebra for all but finitely many $q$.

\begin{theorem}\label{symptheorem}
The decomposition matrix of $\Cc_{n}(q^{-N-1},q)$
with $N>0$ even, for diagrams not fixed by any reflection in $W(B_\infty)$
is given by $(a^{\la}_{\mu}(-N))$, with $a^{\la}_{\mu}(-N)$ as in
Lemma \ref{KLlimit},
for all but finitely many values of $q$.
\end{theorem}

$Proof.$ It has been shown that the algebras $\Cc_n(r,q)$ are
cellular, see \cite{Xi}. This means, in particular, that any
simple $\Cc_n(r,q)$-module is obtained as a quotient of $U_{n,\mu}$
with respect to the annihilator ideal of a certain bilinear form
$(\ ,\ )$. In particular, its dimension is equal to the rank of
the matrix $((b_i,b_j))$ for a certain basis $(b_i)$ of  $U_{n,\la}$.
If $r=q^{-N-1}$, for $N$ positive and even, the entries of this
matrix are Laurent polynomials in $q$, and so are all its minors;
hence its rank is the same for all but finitely many values of $q$.
Using  Lemma \ref{levelrank1}, we know that
$\Cc_n(q^{-N-1},q)\cong \Cc_n(q^{-(2\ell-N-2)-1}, q)$,
for $q$ a primitve $4\ell$-th root of unity. Hence, for $\ell$ sufficiently
large, we have $2\ell-N-2>n$ and we can use Theorem \ref{tiltdual}
to determine the decomposition matrix of these algebras.
It follows from this
and Lemma \ref{KLlimit} that
$a^\la_\mu(-(2\ell-N-2),\ell)=a^\la_\mu(-N)$
for large enough $\ell$. Hence, using  Lemma \ref{levelrank1},
we have for the multiplicity of a simple $\Cc_n(q^{-N-1},q)$ module
labeled by $\la'$ in $U_{n,\mu'}$ that

$$a^\la_\mu(-N,\ell)=a^{\la'}_{\mu'}(-(2\ell-N-2),\ell)=
a^{\la'}_{\mu'}(-N),$$
for sufficiently large $\ell$. As this holds for infinitely
many values of $q$, it must hold for all but finitely many values
of $q$, by the first part of the proof. $\Box$

\medskip

\begin{remark}
1. It is well-known how to calculate the coefficients $a^\la_\mu(-N)$
also for Young diagrams $\la$ which are fixed by a reflection of
$W(B_\infty)$. They can be computed by a slight generalization
of the Kazhdan-Lusztig algorithm, see e.g. \cite{GW2}.

2. Our proof above can be easily adapted to the orthogonal case, i.e. for
the algebras $\Cc_n(q^{N-1},q)$ with $N>0$, assuming an analog
of the result in \cite{Hu} for orthogonal quantum groups. Here the results
would be expressed in terms of parabolic Kazhdan-Lusztig polynomials
for the reflection groups $S_\infty\subset W(D_\infty)$.

3. We expect that we have the same coefficients $a^\la_\mu(-N)$
also for the Brauer algebra $\C_n(-N)=\lim_{q\to 1} \Cc_n(q^{N-1},q)$.
In particular, this should be compatible with the results in \cite{Mt}.

4. It is easy to find a categorical interpretation of our results,
mimicking the theory of tilting modules and fusion categories in the
context of quantum groups: It is well-known how to define a category
$Rep(O(x))$ depending on a paramter $x$, see e.g. \cite{D}. In our
context, this could also be done by applying the idempotent
construction to Brauer algebras (see e.g. \cite{TW2}, \cite{TbW}
for similar constructions in connection with the algebras $\Cc_n(r,q)$).
Then the category $\C(N)$ derived
from the idempotents of the algebras $\C_n(N)$ corresponds
to the category of tilting modules. It contains the category
$Rep(O(N))$ (for $N>0$) and $Rep(Sp(|N|))$ (for $N<0$) as
quotients, similarly as the fusion categories are quotients of
the category of tilting modules of quantum groups.
\end{remark}

\vskip 1cm
\bibliographystyle{plain}

\vskip 1cm

\vskip .3cm
\ni Department of Mathematics, UC San Diego, La Jolla
CA 92093-0112, USA

\ni email:
hwenzl$\char'100$ucsd.edu

\end{document}